\documentclass[preprint,12pt]{elsarticle}

\usepackage{amssymb}
\usepackage{amsmath}

\usepackage{xcolor,amsmath,amsfonts,amsthm,listings,listings-rust,hyperref,caption}
\usepackage[T1]{fontenc}

\hypersetup{
    colorlinks=true,
    citecolor=blue,
    linkcolor=blue,
    filecolor=magenta,      
    urlcolor=cyan,
    pdftitle={Pendleton-Index-Conjecture},
    pdfpagemode=FullScreen
    }

\newcommand{\Z}{\mathbb{Z}}
\newcommand{\R}{\mathbb{R}}

\newtheorem{theorem}{Theorem}[section]
\newtheorem{lemma}[theorem]{Lemma}
\newtheorem{conjecture}[theorem]{Conjecture}

\theoremstyle{definition}
\newtheorem{definition}[theorem]{Definition}

\newtheorem{remark}[theorem]{Remark}

\numberwithin{equation}{section}

\journal{Journal of Number Theory}

\begin{document}

\begin{frontmatter}



\title{Improved Bounds for the Index Conjecture in Zero-Sum Theory}


\author{Andrew Pendleton} 


\begin{abstract}
The Index Conjecture in zero-sum theory states that when $n$ is coprime to $6$ and $k$ equals $4$, every minimal zero-sum sequence of length $k$ modulo $n$ has index $1$. While other values of $(k,n)$ have been studied thoroughly in the last 30 years, it is only recently that the conjecture has been proven for $n>10^{20}$. In this paper, we prove that said upper bound can be reduced to $4.6\cdot10^{13}$, and lower under certain coprimality conditions. Further, we verify the conjecture for $n<1.8\cdot10^6$ through the application of High Performance Computing (HPC).
\end{abstract}



\begin{keyword}
zero-sum theory \sep index conjecture \sep Fourier analysis \sep Euler's totient function


\end{keyword}

\end{frontmatter}



\section{Introduction}
\label{intro}

Let $G$ be an additive cyclic group of order $n$. We say a sequence $(a_1)\dots(a_k)$ over $G$ is \textit{zero-sum} if $\sum_{i}a_i = 0$. It is said \textit{minimal} if it contains no proper, nontrivial subsequence that is itself zero-sum. Given $G\simeq\Z/n\Z$, every $a\in G$ can be thought of as an integer $x\in \{0,1,\dots,n-1\}$. For every integer $y$, we denote $(y)_n$ to be the least nonnegative representative of the congruence class $[y]_n$.

\begin{definition}
    The \textit{index} of a sequence $S=(a_1)\dots(a_k)$ over $G$, written $\text{ind}(S)$, is defined to be \[\min\left\{\dfrac{\sum_{i=1}^{k}(ga_i)_n}{n}:g\in G^{*}\right\}\] where $G^*$ is the set of integers less than and coprime to $n$.
\end{definition}

The study of indices has been ongoing for the past 30 years \cite{lemke_kleitman, geroldinger_early}, and indices were first called such by Chapman et al. in \cite{CHAPMAN1999271}. Since then, connections have been found between the study of indices and the study of integer partitions \cite{Ponomarenko2004MINIMALZS}, atomic monoids and factorization theory \cite{k-atoms-gao, gero-zero-sum}, zero-sum sequences \cite{gero-zero-sum, gryn-book}, Heegard Floer homology \cite{jabuka2011heegaard, GE-dedekind},  Dedekind sums \cite{GE-dedekind, li-vol-1}, and discrepancy estimates \cite{GE2021105410}. The Index Conjecture arises from the study of pairs $(k,n)$ (sequence lengths and moduli) for which every minimal zero-sum sequence is of index $1$. Such $(k,n)$ pairs are called \textit{good} (or, conversely, \textit{bad}).

For $k\le 3$, it is trivial that all $(k,n)$ are good. Further, it is known that for $5\le k\le \frac{n}{2}+1$, all $(k,n)$ are bad. Further still, it was proven in contemporaneous papers by Savchev and Chen \cite{savchev_chen} and Yuan \cite{Yuan2007OnTI} that for $k>\frac{n}{2}+1$, all $(k,n)$ are good.

The case of $k=4$ poses a unique challenge with no such simple characterization. Such pairs with $n$ not coprime to $6$ have been proven bad by Ponomarenko in \cite{Ponomarenko2004MINIMALZS}. However, pairs with $n$ coprime to $6$ have been verified good for $n<1000$ \cite{Ponomarenko2004MINIMALZS}. It is this case to which pertains the Index Conjecture, first posited by Ponomarenko.

\begin{conjecture} (The Index Conjecture)
    For $\gcd(n,6)=1$, every minimal zero-sum sequence $S$ over $G$ of length $4$ has index $1$.
\end{conjecture}

As proven by Shen et al. in \cite{CaixiaShen2014}, to prove the Index Conjecture it suffices to show that the following holds.

\begin{conjecture} (The Index Conjecture \`{a} la Shen et al.)
     For $\gcd(n,6)=1$, let $S = (a_1)(a_2)(a_3)(a_4)$ be a minimal zero-sum sequence over $G$. Suppose $\gcd(n, a_i) = 1$ for all $i$. Then $\textup{ind}(S) = 1$.
\end{conjecture}

The study of these conjectures constitutes an active area of research \cite{GE2021105410,Ge_2018,GE-dedekind,xia-notable,gryn-vishne-notable,shen-xia-notable,li-vol-1,li2013minimal,XIA20134047,zeng2017minimal,CaixiaShen2014}. In particular, Conjecture 1.3 was shown to be true by Ge for $n>10^{20}$ in \cite{GE2021105410}, and this paper seeks to refine the techniques of that work.

\begin{theorem}
    Conjecture 1.3 is true for $n>4.6\cdot 10^{13}$.
\end{theorem}

In Sections 5 and 6, we will show that Conjecture 1.3 also holds when $n<1.8\cdot 10^{6}$ through use of HPC. The additional constraints of this form of the Index Conjecture, in addition to a result of Zeng and Qi \cite{zeng2017minimal} that when $n$ is coprime to $30$, all $(k,n)$ are good, will be essential for the algorithm described in Section 5. 

\section{Notation}
We use $(n,m)$ to denote the greatest common denominator of $n$ and $m$ unless context makes it obvious that we are referring to a pair or interval. We denote the Euler totient function $\phi$, the Ramanujan sum $c_n$, $n$-th Harmonic number $\mathfrak{H}_{n}$, the Riemann Zeta Function $\zeta$, and $g$ to be a generic member of $G^*$.

Further, we define the periodic indicator function \[\chi(t):=\begin{cases}
    1, & \text{if}\;0<\{t\}<1/2\\
    1/2, & \text{if}\;\{t\}=1/2\\
    0, & \text{if}\;\{t\}>1/2
\end{cases}\]
with Fourier coefficients defined to be \[\hat{\chi}(k):=\int_{-1/2}^{1/2}\chi(t)e(-tk)dt\]
which further equals
\[\hat{\chi}(k)=\begin{cases}
    1/2, & \text{if}\;k=0\\
    0, & \text{if}\;k\neq0, 2\mid k\\
    1/(i\pi k), & \text{if}\;2\nmid k
\end{cases}.\]

Define $f$ to be a smoothed version of $\chi$ with Fourier coefficients $\hat{f}$ of finite support $[-H,H]$ such that \[f(x)=\sum_{|h|\le H}\hat{f}(h)e(hx)\] and $\hat{f}(0)=1/2$, $\hat{f}(k)=0$ for $k\neq0, 2\mid k$, and $|\hat{f}(k)|\le 1/(\pi(|k|)$ for odd $k$. Furthermore, either $\hat{f}(k)\in\R$ or $\hat{f}(k)\in i\R$.\footnote{I would encourage the reader to refer to the notation section of \cite{GE2021105410} if they are curious about the derivation of these properties and their Fourier analytic rationale.}

\section{Lemmas}
\begin{lemma}
    Let $(n,b) = (n, 6) = 1$ and $H$ be a fixed positive integer. Let $h_2$ be an odd integer and define $h_2^*$ (depending on $b, n, H$) be the only possible integer in $[-H^2, H^2]$ such that $(bh_2 + h_2^*, n) >\sqrt{2H^2n}$ (if such an integer exists). Let $S$ (also depending on $b, n, H$) denote the set $\{(h_2, h_2^*) : |h_2| \le H, 2 \nmid h_2, |h_2^*| \le H, 2 \nmid h_2^*\}$. Assume that none of $b \pm 1$, $3b \pm 1$, $b \pm 3$ is congruent to 0 modulo n. Then for
    \[S_b^*:=\sum_{S}\hat{f}(h_2)\hat{f}(h_2^{*})c_n(bh_2+h_2^*)\] we have \[S_b^*\le 0.079021\cdot\phi(n).\]
\end{lemma}

\begin{proof}
    Note that aside from the final coefficient of $\phi(n)$, this is similar to Lemma 6 in \cite{GE2021105410}. Further note that in said Lemma, Case 2, Subcase 2 bounds the remaining cases and subcases, with Case 1 providing a $\phi(n)$ coefficient of around $0.047$ and Case 2, Subcase 1 providing a $\phi(n)$ coefficient of $0.05$.
    
    We can reduce the bound in Subcase 2 by observing $$|c_n(3^m b+(3^m)^*)|\leq \phi(n)/4,\text{ for all }m\in\Z^{\ge 0}.$$ Therefore, we have $$|\sum_{S}\hat{f}(h_2)\hat{f}(h_2^*)c_n(bh+h_2^*)|$$ $$\le2\sum_{m}|\hat{f}(3^m)\hat{f}((3^m)^*)
    c_n(3^mb+(3^m)^*)|+2\sum_{h_2\neq 3^m,2\nmid h_2}|\hat{f}(h_2)\hat{f}(h_2^*)|\phi(n)$$ $$\le 2\cdot \dfrac{\phi(n)}{4}\sum_{m}\dfrac{1}{3^{2m}\pi^2}+2\cdot\dfrac{1}{\pi^2}(\dfrac{\pi^2}{8}-\sum_{m}\frac{1}{3^{2m}})\cdot\phi(n)$$$$=(\dfrac{9}{16\pi^2}+\dfrac{1}{4}-\dfrac{9}{4\pi^2})\cdot\phi(n)\le 0.079021\cdot\phi(n).$$
\end{proof}

\begin{lemma}
    Let $(b,n)=(n,6)=1$ and $H$ be a fixed positive integer. For any $h_2\in[-H,H]$ there exists at most one $\tilde{h_2}$ (depending on $b,n,H$) in the same domain such that $(bh_2+\tilde{h_2},n)>\sqrt{2Hn}$. Furthermore, for any $\tilde{h_2}\in[-H,H]$ there is at most one $h_2$ in the same domain such that $(bh_2+\tilde{h_2},n)>\sqrt{2Hn}$.
    \end{lemma} 

\begin{proof}
    Define $\mathcal{H}\subseteq[-H,H]$ to be the integers $x$ such that there exists at least one $y\in[-H,H]$ such that \[(bx+y,n)>\sqrt{2Hn}\] and $\tilde{\mathcal{H}}\subseteq[-H,H]$  to be the integers $y$ such that there exists at least one $x\in[-H,H]$ satisfying the same inequality. 
    
    By the first statement of Lemma 5 of \cite{GE2021105410}, there exists a well-defined map \[\alpha:\mathcal{H}\to[-H,H]\] \[h_2\mapsto \tilde{h_2}\] such that $(bh_2+\tilde{h_2},n)>\sqrt{2Hn}$. To prove the lemma, it suffices to show $\alpha$ is injective.

    By the second statement of the same Lemma 5, there exists a well-defined map \[\beta:\tilde{\mathcal{H}}\mapsto[-H,H]\] \[\tilde{h_2}\mapsto h_2\] such that $(bh_2+\tilde{h_2},n)>\sqrt{2Hn}$.

    Suppose, for the sake of contradiction, we have $h_2,h_2'\in\mathcal{H}$ not equal such that $$\alpha(h_2)=\alpha(h_2').$$ By definition we have \begin{equation}
    (bh_2+\alpha(h_2),n)>\sqrt{2Hn}
    \end{equation} and \begin{equation}
    (bh_2'+\alpha(h_2),n)>\sqrt{2Hn}.
    \end{equation} By (3.1), $\alpha(h_2)\in\tilde{\mathcal{H}}$. So, consider $\beta(\alpha(h_2))$. By (3.1) we have that $\beta(\alpha(h_2))=h_2$ and by (3.2) we have $\beta(\alpha(h_2))=h_2'$. However, $h_2\neq h_2'$ by assumption, a contradiction with $\beta$ well-defined. Therefore, no such $h_2,h_2'$ can exist, and consequently, $\alpha$ is injective. Indeed, the corestriction of $\alpha$ to $\tilde{\mathcal{H}}$ is a bijection.
\end{proof}

\begin{lemma}
    Let $(n, b) = (n, 6) = 1$ and $H$ be a fixed positive integer. Let $h_2$ be an odd integer and define $\tilde{h_2}$ (depending on $b, n, H$) to be the only possible integer in $[-H,H]$ such that $(bh_2+\tilde{h_2},n)>\sqrt{2Hn}$. If such an integer exists, it is unique by Lemma 3.2. Let $T$ (also depending on $b, n, H$) denote the set $\{(h_2, \tilde{h_2}) : |h_2| \le H, 2 \nmid h_2, |\tilde{h_2}| \le H, 2 \nmid \tilde{h_2}\}$. Then for \[\tilde{T_b}:=\sum_{T\setminus S}\hat{f}(h_2)\hat{f}(\tilde{h_2})c_n(bh_2+\tilde{h_2}),\] where $S$ is defined as in Lemma 3.1, we have\[|\tilde{T_b}|\le \frac{\sqrt{2H^2n}}{4}.\] Further, if any of $3b\pm 1$ or $b\pm 3$ is congruent to $0$ modulo $n$, then \[|\tilde{T_b}|\le \frac{\sqrt{2H^2n}}{12}.\]
\end{lemma}

\begin{proof}
    By Lemma 3.1, it is certainly true that $(bh_2+\tilde{h_2},n)<\sqrt{2H^2n}$ for any $(h_2,\tilde{h_2})\in T\setminus S$. So, by the Cauchy-Schwartz inequality we have \[|\tilde{T_b}| \le \sqrt{2H^2n}\left(\sum_{0<|h_2| \le H}|\hat{f}(h_2)||\hat{f}(\tilde{h_2})|\right) \le \sqrt{2H^2n}\left(\sum_{h_2}|\hat{f}(h_2)|^2\right)\]\[ \le \frac{2\cdot\sqrt{2H^2n}}{\pi^2}\left(\sum_{k=1}^\infty\frac{1}{(2k-1)^2}\right)=\frac{\sqrt{2H^2n}}{4}.\] The last equality holds as a consequence of $$\zeta(2)=\sum_{k=1}^\infty\frac{1}{k^2}=\sum_{k}\frac{1}{(2k)^2}+\sum_{k}\frac{1}{(2k-1)^2}=\frac{1}{4}\cdot\zeta(2)+\sum_{k}\frac{1}{(2k-1)^2}$$ $$\implies \sum_{k}\frac{1}{(2k-1)^2}=\frac{3}{4}\cdot\zeta(2)=\frac{\pi^2}{8}.$$ As for the second case, without loss of generality, assume $b+3\equiv0\pmod{n}$, then \[|\tilde{T_b}|\le\dfrac{1}{3}\cdot\sqrt{2H^2n}\cdot\dfrac{2}{\pi^2}\cdot\left(\sum_{k=1}^\infty\frac{1}{(2k-1)^2}\right)=\dfrac{\sqrt{2H^2n}}{12}.\]
\end{proof}

\begin{lemma}
    Let $(n, b) = (n, 6) = 1$ and $H$ be a fixed positive integer. Let $h_2,h_3\in[-H,H]$ such that $(bh_2+h_3,n)<\sqrt{2Hn}$. Then for
    \[R_b:=\sum_{h_2}\sum_{h_3}\hat{f}(h_2)\hat{f}(h_3)c_n(bh_2+h_3)\]we have\[|R_b|\le \sqrt{2Hn}\left(\frac{2}{\pi}(\mathfrak{H}_{2\theta}-\frac{1}{2}\mathfrak{H}_\theta)\right)^2\]
\end{lemma}

\begin{proof}
    Define $\theta:=\lceil H/2\rceil$. For $R_b$ we have \[|R_b|\le\sqrt{2Hn}\left(\sum_{h_2}\sum_{h_3}|\hat{f}(h_2)||\hat{f}(h_3)|\right)\le\sqrt{2Hn}\left(\sum_{h_2}|\hat{f}(h_2)|\right)^2\]\[\le\sqrt{2Hn}\left(\dfrac{2}{\pi}\sum_{1\le h_2\le H}\frac{1}{h_2}\right)^2\le\sqrt{2Hn}\left(\frac{2}{\pi}\sum_{k=1}^{\theta}\frac{1}{2k-1}\right)^2\]\[=\sqrt{2Hn}\left(\frac{2}{\pi}(\mathfrak{H}_{2\theta}-\frac{1}{2}\mathfrak{H}_\theta)\right)^2.\]
\end{proof}

\begin{lemma}
    Let $P_n:=\{p_1,p_2,\dots,p_k\}$ be a finite list of distinct primes coprime to and less than $n$ and let $p(n):=\prod_{P_n}p_i$. Then, \[\phi(n)>\dfrac{n\cdot p(n)}{\phi(p(n))\left( e(\gamma)\log\log(n\cdot p(n))+\dfrac{2.50637}{\log\log(n\cdot p(n))}\right)},\] where $\gamma$ is Euler's constant.
\end{lemma}

\begin{proof}
    As the Euler totient is known multiplicative for coprime factors, \begin{equation}
        \phi(n\cdot p(n))=\phi(n)\phi(p(n)).
    \end{equation} Using the lower bound for $\phi(n\cdot p(n))$ found by Rosser and Schoenfeld \cite{ross_scho}, \begin{equation}
        \phi(n\cdot p(n))>\dfrac{n\cdot p(n)}{e(\gamma)\log\log(n\cdot p(n))+\dfrac{2.50637}{\log\log(n\cdot p(n))}}.
    \end{equation}
    The premise follows by (3.3) and (3.4).
\end{proof}
Given $(6,n)=1$, we can specifically say that \[\phi(n)>\dfrac{3n}{e(\gamma)\cdot\log\log(6n)+\dfrac{2.50637}{\log\log(6n)}}.\] We prove this general formulation to provide smaller bounds under additional coprimality conditions of $n$ in Remark 4.1.
\begin{lemma}
    For integers $a,b$ comprime to $n$, define \[S_0:=\sum_g\chi(g/n)\chi(ag/n)\chi(bg/n)\] and \[S_1:=\sum_g f(g/n)f(ag/n)f(bg/n).\] For $H>1000$ we have \[|S_0-S_1|\le \frac{13.02}{H}\cdot\phi(n)+20.02\sqrt{2Hn}+7H.\]
\end{lemma}

\begin{proof}\renewcommand{\qedsymbol}{}
    See Lemma 8 of \cite{GE2021105410}.
\end{proof}

In fact, given our future choice of $H$, we may assume that $H>5000$, but this yields negligible improvement to our final bound.

\begin{lemma}
    Let $(n,6)=1$ and $g\in G^*$. For a minimal zero-sum sequence modulo $n$ with index $2$, exactly two of $(g)_n,(ga)_n,(gb)_n,(gc)_n$ lie in the interval $(0,n/2)$.
\end{lemma}

\begin{proof}\renewcommand{\qedsymbol}{}
See Remark 2.1 of Li and Peng \cite{li2013minimal}.    
\end{proof}

Now observe that the index of a minimal zero-sum 4-sequence is either $1$ or $2$ because if a given $g$ makes \[\dfrac{gi+ga+gb+gc}{n}=3\] then \[\dfrac{(n-g)i+(n-g)a+(n-g)b+(n-g)c}{n}=1.\] Further, this implies \[c=2n-a-b-1.\] This lemma informs the general strategy of the final proof; when the modulus $n$ is sufficiently large, more than two of $(g)_n,(ga)_n,(gb)_n,(gc)_n$ lie in the interval $(0,n/2)$.

\begin{lemma}
    Let $S=(1)(a)(b)(c)$, where $c = 2n-1-a-b$, be a minimal zero-sum sequence modulo $n$, whose index is $2$. If $n>100$, then at least two of the pairs $(1,a),(1,b),(a,b)$  do not satisfy any of the linear relations $x \pm 3y \equiv 0 \pmod{n}$ or $3x\pm y \equiv 0 \pmod{n}$ for $(x,y)$.
\end{lemma}

\begin{proof}
By Lemma 9 in \cite{GE2021105410} we know that at least one of the pairs $(1,a)$, $(1,b)$, $(a,b)$ does not satisfy any linear relations \begin{equation}
    x \pm 3y \equiv 0 \pmod{n}
\end{equation} or \begin{equation}
    3x \pm y \equiv 0 \pmod{n}
\end{equation} for $(x,y)$. What remains to be shown is that two of the pairs $(1,a)$, $(1,b)$, $(a,b)$ cannot satisfy any of (3.5) or (3.6) simultaneously for $(x,y)$. Proving this amounts to demonstrating a contradiction in each of the $\binom{3}{2}\cdot 4^2=48$ possible cases. For the majority of cases this is relatively straightforward, if tedious work. Generally speaking, the information provided suffices to fix $S$. However, there are four cases that prove more subtle than the rest. These cases are when \[1 \pm 3a\equiv 0 \pmod{n}~\text{and}~a-3b\equiv 0 \pmod{n}\] or when \[1 \pm 3b\equiv 0 \pmod{n}~\text{and}~3a+b\equiv 0 \pmod{n}.\] In the interest of brevity, only the case when \begin{equation}
    1-3b\equiv 0 \pmod{n}
\end{equation} and \begin{equation}
    3a+b\equiv 0 \pmod{n}
\end{equation} will be detailed here. The other cases remain largely similar. First, from (3.7), observe that \[b=\dfrac{1-kn}{3}\] for some $k\in\Z$. From Lemma 3.7 and Theorem 1.1 of \cite{Ge_2018} we know that \[n/2<b<n\] for a minimal zero-sum sequence of index 2. Therefore, given $n>100$, $k=-2$ and \begin{equation}
    b=\dfrac{2n+1}{3}.
\end{equation} (Recall that the Index Conjecture has been computationally verified for $n<1000$ \cite{Ponomarenko2004MINIMALZS}.) Substituting (3.9) into (3.8) and simplifying we have that \[a=\dfrac{ln-1}{9}\] for some $l\in\Z$. Again, from Lemma 3.7 and Theorem 1.1 of \cite{Ge_2018} we know that \[1<a<n/2\] for a minimal zero-sum sequence of index 2. Therefore, yet again noting $n>100$, $l=4$ and \begin{equation}
    a=\dfrac{4n-1}{9}.
\end{equation} Given $c=2n-1-a-b$, we can now calculate $c$ from (3.9) and (3.10) explicitly as \begin{equation}
    c=\dfrac{8n-11}{9}.
\end{equation} Since $a,b,c\in\Z^+$ and $(n,6)=1$, we can derive a system of linear relations for $n$, \[\begin{cases}
    4n-1\equiv 0 \pmod{9}\\
    2n-1\equiv 0 \pmod{3}\\
    8n-11\equiv 0 \pmod{9}\\
    n\equiv 1 \pmod{2}\\
    n\equiv 0 \pmod{5}
    \end{cases}\] which simplifies to \[\begin{cases}
        n\equiv 7\pmod{9}\\
        n\equiv 0\pmod{5}\\
        n\equiv 1\pmod{2}
    \end{cases}.\] From this, we can derive a formula for $n$ by elementary application of the Chinese remainder theorem, \begin{equation}
        n=90m+25
    \end{equation} for some $m\in\Z^{\geq 0}$. Combining (3.9), (3.10), (3.11), and (3.12) we have that \[S=(1)(40m+11)(60m+17)(80m+21)\] for some $m\in\Z^{\geq 0}$. But, such a sequence cannot have index 2 because for $g=17$, \[\dfrac{(g)_n+(ga)_n+(gb)_n+(gc)_n}{n}=\dfrac{17+50m+12+30m+14+10m-18}{90m+25}=1.\] There exist such special $g$ values for the other subtle cases, as well. In particular, $g=4$ suffices when \[1 + 3a\equiv 0 \pmod{n}~\text{and}~a-3b\equiv 0 \pmod{n},\] $g=3$ suffices when \[1 - 3a\equiv 0 \pmod{n}~\text{and}~a-3b\equiv 0 \pmod{n},\] and $g=17$ suffices when \[1+3b\equiv 0 \pmod{n}~\text{and}~3a+b\equiv 0 \pmod{n}.\] As every case results in a contradiction, at least two of the pairs $(1,a),(1,b),(a,b)$ do not satisfy any of the linear relations $x \pm 3y \equiv 0 \pmod{n}$ or $3x\pm y \equiv 0 \pmod{n}$ for $(x,y)$. 
\end{proof}

\section{Proof of Theorem 1.4}
By Lemma 3.7, a sequence has index 1 if there exists a $g\in G^*$ such that $3$ of $\frac{g}{n}$, $\frac{(ga)_n}{n}$, $\frac{(gb)_n}{n}$, $\frac{(gc)_n}{n}$ are in $(0,1/2)$. Recall \[S_0=\sum_{g}\chi(g/n)\chi(ag/n)\chi(bg/n)\] and \[S_1=\sum_{g}f(g/n)f(ag/n)f(bg/n).\] To prove the theorem, it suffices to show that $S_0>0$. Specifically, we will show that $ S_1-(S_1-S_0)>0$.

We begin by finding a lower bound for $S_1$. First note that, as $f$ is real-valued, $S_1$ is real and so $S_1=\mathfrak{R}S_1$. Rewriting $S_1$ in terms on $\hat{f}$ along support $H$, we have \[S_1=\sum_{h_1}\sum_{h_2}\sum_{h_3}\hat{f}(h_1)\hat{f}(h_2)\hat{f}(h_3)\sum_{g}e\left(\dfrac{g}{n}(ah_1+bh_2+h_3)\right).\] Suppose, that exactly two of $h_1, h_2, h_3=0$. Without loss of generality, assume $\hat{f}(h_1)=\hat{f}(h_2)=1/2$ and $\hat{f}(h_3)=0$ or is in $i\R$. As the Ramanujan sum is always real, \[\hat{f}(h_1)\hat{f}(h_2)\hat{f}(h_3)\sum_ge\left(\dfrac{g}{n}(ah_1+bh_2+h_3)\right)=0\] or is in $i\R$. Suppose similarly, that none of $h_1, h_2, h_3=0$. For every $h\in\{h_1,h_2,h_3\}$, $\hat{f}(h)=0$ or is in $i\R$. In either case, $S_1=\mathfrak{R}S_1=0$, so we need only consider when either all of $h_1, h_2, h_3=0$ or exactly one of the $h$ values is $0$, as these are the only terms that survive.

Further, \begin{equation}
    S_1=\phi(n)\cdot(\hat{f}(0))^3+\hat{f}(0)\cdot\left(\sum_{h_2}\sum_{h_3}+\sum_{h_3}\sum_{h_1}+\sum_{h_1}\sum_{h_2}\right),
\end{equation} where \[\sum_{h_2}\sum_{h_3}=\sum_{0<|h_2|\le H}\sum_{0<|h_3|\le H}\hat{f}(h_2)\hat{f}(h_3)\sum_{g}e(\frac{g}{n}(bh_2+h_3))\] and the other double sums are defined similarly. We break this sum into three pieces \[\sum_{h_2}\sum_{h_3}=S_b^*+\tilde{T_b}+R_b.\] The part $S^*$ is the sum over \[S=\{(h_2, h_2^*) : |h_2| \le H, 2 \nmid h_2, |h_2^*| \le H, 2 \nmid h_2^*,\}\] as defined in Lemma 3.1. The part $\tilde{T_b}$ is the sum over \[T\setminus S=\{(h_2, \tilde{h_2}) : |h_2| \le H, 2 \nmid h_2, |\tilde{h_2}| \le H, 2 \nmid \tilde{h_2}\}\setminus S,\] as defined in Lemma 3.2. The final part $R_b$ is the double sum over the the remaining $(h_2,h_3)$ pairs.

By Lemma 3.8, at least two of the pairs $(1, a),(1, b),(a, b)$ do not satisfy any linear relations \begin{equation}
    x \pm 3y \equiv 0 \pmod{n}
\end{equation} or \begin{equation}
    3x \pm y \equiv 0 \pmod{n}
\end{equation} for $(x, y)$. Without loss of generality, take $(1,a)$ and $(1,b)$ to be the two such pairs. Therefore, $|S_a^*|,|S_b^*|$ are bounded by Lemma 3.1 and $|\tilde{T_a}|,|\tilde{T_b}|$ by the first case of Lemma 3.3.

For the remaining sums, there are two cases. If $(a,b)$ satisfies neither (4.2) nor (4.3), then similarly to $|S_a^*|$ and $|\tilde{T_a}|$ we have by Lemma 3.1 \[
    |S_{ab^{-1}}^*|\le0.079021\cdot\phi(n). 
\] and by the first case of Lemma 3.3 \[
    |\tilde{T_{ab^{-1}}}|\le\dfrac{1}{4}\cdot\sqrt{2H^2n}. 
\] If $(a,b)$ satisfies any of (4.2) or (4.3), then by Lemma 7 in \cite{GE2021105410} we have \begin{equation}
    |S_{ab^{-1}}^*|\le\frac{1}{12}\cdot\phi(n).
\end{equation} and by Lemma 3.3, \begin{equation}
    |\tilde{T_{ab^{-1}}}|\le\dfrac{1}{12}\cdot\sqrt{2H^2n}.
\end{equation}

Assume we have the second case. By the bounds provided in Lemma 3.1, the first case of Lemma 3.3, and Lemma 3.4 we have \begin{equation}
\left|\sum_{h_2}\sum_{h_3}\right|\le 0.079021\cdot\phi(n)+\frac{\sqrt{2H^2n}}{4}+\sqrt{2Hn}\left(\frac{2}{\pi}(\mathfrak{H}_{2\theta}-\frac{1}{2}\mathfrak{H}_\theta)\right)^2
\end{equation} where $\theta=\lceil H/2\rceil.$ Note that we can bound the final double sum over $h_3$ and $h_1$ in the exact same manner. We proceed similarly for the double sum over $h_1$ and $h_2$. By (4.4), (4.5), and Lemma 3.4 we have \begin{equation}
    \left|\sum_{h_1}\sum_{h_2}\right|\le\dfrac{1}{12}\cdot\phi(n)+\frac{\sqrt{2H^2n}}{12}+\sqrt{2Hn}\left(\frac{2}{\pi}(\mathfrak{H}_{2\theta}-\frac{1}{2}\mathfrak{H}_\theta)\right)^2.
\end{equation} As the Ramanujan sum may be negative, we assume each $\sum_{h_2}\sum_{h_3}$ is negative. So, by (4.1), (4.6), and (4.7) we have \begin{equation}
    S_1\ge c_0\phi(n)-\frac{1}{2}(\dfrac{1}{4}+\dfrac{1}{4}+\dfrac{1}{12})\cdot\sqrt{2H^2n}-\dfrac{3}{2}\left(\sqrt{2Hn}\left(\frac{2}{\pi}(\mathfrak{H}_{2\theta}-\frac{1}{2}\mathfrak{H}_\theta)\right)^2\right),
\end{equation} where \[c_0=\frac{1}{8}-\frac{1}{2}(0.079021+0.079021+\frac{1}{12}).\]

Thus, by (4.8) and Lemma 3.6 we have
    \[S_0\ge c_{1}\phi\left(n\right)-\frac{7}{24} \sqrt{2H^2n}-\left(\frac{3}{2}\left(\frac{2}{\pi}(\mathfrak{H}_{2\theta}-\frac{1}{2}\mathfrak{H}_\theta)\right)^{2}+20.02\right)\sqrt{2Hn}-7H,\]
where \[c_1=c_0-\frac{13.02}{H}.\] We can optimize our choice of $H$ by minimizing the ratio $$\dfrac{H}{c_1}=\dfrac{H}{c_0-\frac{13.02}{H}}$$ such that both the numerator and denominator are positive. Doing so via standard analytic methods gives $H\approx 7000$, which we take as our $H$. Further, take $c_{0}-\frac{13.02}{H}\approx 0.0024523$ to be our $c_1$. Note that such a choice for $H$ suffices for Lemma 3.6.

Using the lower bound for $\phi(n)$ in Lemma 3.5 with $P_n=\{2,3\}$, one can calculate explicitly that it suffices to take $n>4.6\cdot 10^{13}.$ One can now calculate and see that the previous assumption regarding (4.4) and (4.5) is indeed sensible; the other case will result in less restrictive bounds upon $n$. \qed

\begin{remark}
    Applying the general form of Lemma 3.5, the following improved upper bounds can be computed for $n$ satisfying additional coprimality conditions. Again we take $H=7000$ and $c_1=0.0024523$.
    
    In fact, if we assume $n$ is only divisible by powers of $5$, then we can infer $\phi(n)=4n/5.$ Performing the same computation again, we see the following upper bounds are bounded below by $\approx 1.4\cdot 10^{13}$. Thus, if we wish to further lower the upper bound for $n$ where $P_n=\{2,3\}$, there is relatively little to be gained by computing closer approximations of $\phi$. The case of nearly prime $n$ has been studied extensively. Specifically, the Index Conjecture has recently been proven true for $n$ with two prime factors \cite{XIA20134047}.
\end{remark}

\begin{table}[!ht]
    \centering
    \begin{tabular}{cc}
        $P_n$ & Upper Bound for $n$ \\\hline
        $\{2,3\}$ & $4.6\cdot 10^{13}$ \\
        ---\textquotedbl---$\;\cup\{7\}$ & $3.4\cdot 10^{13}$\\
        ---\textquotedbl---$\;\cup\{11\}$ & $2.9\cdot 10^{13}$\\
        ---\textquotedbl---$\;\cup\{13\}$ & $2.6\cdot 10^{13}$\\
        ---\textquotedbl---$\;\cup\{17\}$ & $2.3\cdot 10^{13}$\\
        ---\textquotedbl---$\;\cup\{19\}$ & $2.2\cdot 10^{13}$\\
    \end{tabular}
    \caption*{Additional Upper Bounds for Special $n$}
    \label{tab:my_label}
\end{table}

\section{Computing Strategy}
The first to take a computational approach to the Index Conjecture was Ponomarenko, who verified that it is indeed true for $n<1000$ in 2004 \cite{Ponomarenko2004MINIMALZS}. Given the development of more powerful and accessible computational resources in the past $20$ years, we can push this lower bound much further. We make use of William \& Mary's high performance computing (HPC) resources (specifically 8-16 nodes of the \href{https://www.wm.edu/offices/it/services/researchcomputing/hw/nodes/kuro/}{Kuro cluster}) to verify that Conjecture 1.3 holds for $n<1.8\cdot10^6$.

First, we generate a list of $n$ values which we wish to check (where $(n,6)=1$ and $5\mid n$, recall that \cite{zeng2017minimal} gives us the latter divisibility condition). By an elementary application of the Principle of Inclusion and Exclusion, it can be seen that there are $c(M)-c(N)$ such $n$ values between positive integers $N$ and $M$ where \[c(N):=-\left(\lfloor\dfrac{N}{2}\rfloor+\lfloor\dfrac{N}{3}\rfloor-\lfloor\dfrac{N}{5}\rfloor\right)+\left(\lfloor\dfrac{N}{6}\rfloor+\lfloor\dfrac{4}{5}\cdot\lfloor\dfrac{N}{2}\rfloor\rfloor+\lfloor\dfrac{4}{5}\cdot\lfloor\dfrac{N}{3}\rfloor\rfloor\right)-\lfloor\dfrac{4}{5}\cdot\lfloor\dfrac{N}{6}\rfloor\rfloor\]\[\approx\lfloor\dfrac{N}{15}\rfloor.\] We use \href{https://www.lustre.org/}{Lustre} to locally store files containing all $n$ values to be checked, those being checked, and those which have been checked. Then, each thread on each node (16 nodes equates to roughly 1024 threads working simultaneously) pulls an $n$ value from the list and checks it---after which it marks that $n$ value as complete, and grabs another.

``Checking" an $n$ value amounts to iterating over the possible minimal zero-sum 4-sequences $$i,a,b,c$$ for the given modulus $n$ and calculating their indices. Therefore, we need only check sequences that could have an index different from $1$.

First, note that we can multiply the sequence by $i^{-1}$ to fix the first element as $1$. Again, note that a minimal zero-sum sequence of length $4$ is either $1$ or $2$. We choose $1<a<b<c<n-1$, ensuring the sequence is minimal. (As $1+a<1+b<1+c<n$.) By Lemma 3.7, we may assume \begin{equation}
    2\le a< n/2 < b,
\end{equation} and by Theorem 1.1 of \cite{Ge_2018}, we may further assume \begin{equation}
     b+1<c<n-2.
\end{equation} Therefore, we have ranges for possible $a$ and $b$ values, namely \begin{equation}
    n+2-a \le b\le n-\dfrac{3}{2}-\dfrac{a}{2}
\end{equation} and \begin{equation}
    7\le a< \dfrac{n}{2}.
\end{equation}

Accordingly, in lines 12 to 22 in the following code, we define a function that computes the index of a given sequence given its modulus and a vector of all numbers less than and coprime to the modulus.

In lines 24 to 71, we define a function that takes in a modulus, computes a vector of all numbers less than and coprime to the modulus, and then using the conditions described in (5.1), (5.2), and (5.4), filters that vector further into vectors of possible $a$ and $b$ values. We then iterate over the possible $a$ values---ensuring that we only check $b$ values satisfying (5.3)---to check that the sequence $(1)(a)(b)(c)$ is indeed zero-sum, and call the previously defined function to calculate the index. If it is not $1$, it breaks.

Lastly, lines 73 to 203 have been omitted. These lines moderate the interaction between the process instance and the lists of stored $n$ values, specifically written to avoid a race condition upon retrieval. Whereupon, the $n$ values are distributed to threads using the \href{https://docs.rs/rayon/latest/rayon/}{Rayon} library.

We choose to write this algorithm in Rust for its speed and support for multithreading, which further contributes to the speed of computation, as compared to another language like Python that is not a compiled language and is less amenable to multithreading.

\section{Code}

\begin{lstlisting}[basicstyle=\footnotesize, language=Rust, numbers=left, breaklines=true]
use hashbrown::HashSet;
use num::Integer;
use rayon::prelude::*;

use std::collections::VecDeque;
use std::fs;
use std::io::ErrorKind;
use std::path::Path;
use std::thread;
use std::time::Duration;

fn w_index(s: [i64; 4], n: i64, coprimes: &[i64]) -> i64 {
    for g in coprimes {
        let sum: i64 = s.iter().map(|i| (i * g) % n).sum();
        let sum_div = sum / n;
        if sum_div != 2 {
            return 1;
        }
    }

    2
}

fn big_check(n: i64) {
    let half_n = n / 2;
    let coprimes: Vec<i64> = (1..n).into_par_iter().filter(|&i| i.gcd(&n) == 1).collect();

    let coprimes_a: Vec<i64> = (&coprimes)
        .into_par_iter()
        .filter_map(|i| {
            if *i >= 7 && *i < half_n {
                Some(*i)
            } else {
                None
            }
        })
        .collect();

    let coprimes_b: Vec<i64> = (&coprimes)
        .into_par_iter()
        .filter_map(|i| {
            if *i > half_n && *i < n - 2 {
                Some(*i)
            } else {
                None
            }
        })
        .collect();

    let coprime_set: HashSet<&i64> = HashSet::from_iter(coprimes.iter());

    coprimes_a.into_par_iter().for_each(|a| {
        for &b in coprimes_b.iter() {
            if b >= n + 2 - a && b <= n - (3 / 2) - (a / 2) {
                let c = 2 * n - a - b - 1;

                if coprime_set.contains(&c) {
                    let s = [1, a, b, c];
                    let sum: i64 = s.iter().sum();

                    if sum % n == 0 {
                        if w_index(s, n, &*coprimes) != 1 {
                            println!("error at: {} for: {:?}", n, s);
                            return;
                        }
                    }
                }
            }
        }
    })
}
\end{lstlisting}

\section*{Acknowledgements}
The author would like to thank Professor Fan Ge for suggesting the Index Conjecture as the topic of this honors thesis and for spending his Thursday afternoons providing feedback, resources, and revisions. He also acknowledges Thomas Schollenberger of the Rochester Institute of Technology for helping rewrite the author's glacial Python code into Rust and assisting with the operation of the HPC file system and terminal. Finally, the author is grateful to those who attended the original defense of this paper, including but not limited to: Professor Pierre Clare and Professor Matthew Schueller. He would also like to thank the referees for their contributions.

The author acknowledges William \& Mary Research Computing for providing computational resources and/or technical support that have contributed to the results reported within this paper. URL: \url{https://www.wm.edu/it/rc}
\bibliographystyle{elsarticle-num} 
\bibliography{pendletonIC_citations}






\end{document}